# On one method of boundary value problem regularization by passage to the limit

**Vladimir Gotsulenko, Lyudmila Gaponova**
"Strategy" Institute for Entrepreneurship, Zhovti Vody, Ukraine

**ABSTRACT.** For one class of boundary value problem depending on small parameter for which numerical methods for their solution are actually inapplicable, procedure of limiting problem acquisition which is much easier and which solution as much as close to the initial solution is described.

For the objective description of many physical, social and economic, biological and other processes their mathematical models are created. Many mathematical models consist of, as a rule, nonlinear differential equation or their system and some unambiguity conditions (i.e. initial and boundary conditions). As only for very few classes of the differential equations analytical methods of their integration are developed, for obtaining of result actually always it is necessary to use the numerical methods, for example widespread net method, by replacing derivatives by corresponding finite differences. As a result all is reduced to the solution of system of the algebraic equations. However there are problems [КГ05] for which the quantity of the equations in system is so great, even with the requirement of small accuracy of the numerical solution, that obviously any computer memory will not be enough for their solution. To the solution, to be exact to regularization of one class of such problems the given work is devoted. Let there is a boundary value problem:

$$Lu_\varepsilon = f \text{ On } \Omega_\varepsilon \quad [1]$$

$$u_\varepsilon \big|_{\partial\Omega_\varepsilon} = a_\varepsilon \quad [2]$$

where $\Omega_\varepsilon$ is $N$-connected punched set and it is formed of some simply connected open set $\Omega \subset R^q$ by removal from it periodically located closed subsets $T_k^\varepsilon$, i.e. $\Omega_\varepsilon = \Omega - \bigcup_{k=1}^{N} T_k^\varepsilon$ (figure 1), $L$ - some differential operator $u_\varepsilon \in Y(\Omega_\varepsilon)$, $Y(\Omega_\varepsilon)$ - separable reflexive Banach or Hilbert space.

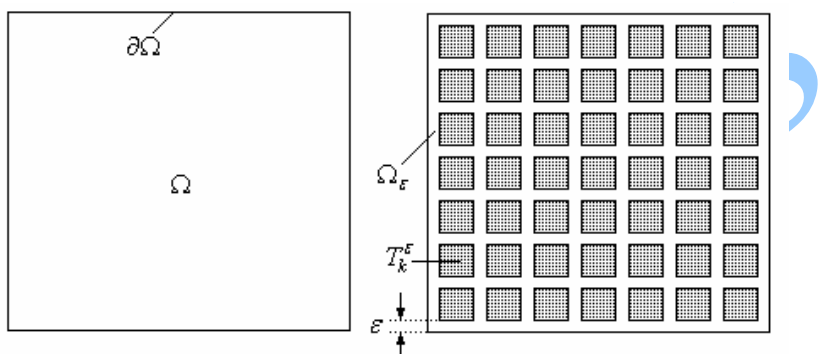

***Fig. 1*** *Punching of set $\Omega$ by subsets $T_k^\varepsilon$*

Let's note, that a boundary value problem [1] - [2] with $\varepsilon \to 0$ is not possible to solve by any numerical method since the system of the algebraic equations is obtained where quantity $N \approx O\left(\frac{1}{\varepsilon}\right) \to \infty$. Understanding that function $u_\varepsilon \in Y(\Omega_\varepsilon)$ is the weak solution of the differential equation [1], it is possible to receive a following aprioristic estimation being based on integrated Green type formula:

$$\|u_\varepsilon\|_{Y(\Omega_\varepsilon)} \leq \gamma_o \|f\| + \gamma_1 \|a_\varepsilon\| \leq const \qquad [3]$$

Further, using the advantage of the known theorem that the limited set in separable reflexive Banach space is weak precompact, we shall receive, that

$$\tilde{u}_\varepsilon \to \tilde{u}_o \text{ weak in } Y(\Omega) \qquad [4]$$

where function $\tilde{u}_\varepsilon = \begin{cases} u_\varepsilon, \Omega_\varepsilon \\ \xi_\varepsilon, \bigcup_{k=1}^{N} T_k^\varepsilon \end{cases}$, and functions $\xi_k$ get out of a condition that the investment $\tilde{u}_\varepsilon \in Y(\Omega)$ took place.

Carrying out passage to the limit with $\varepsilon \to 0$ in a task [1] - [2], as a result we receive a boundary value problem:

$$M\tilde{u}_o = f \text{ On } \Omega \qquad [5]$$
$$\tilde{u}_o|_{\partial\Omega} = b \qquad [6]$$

where $M$ - the differential operator, generally speaking, differing the operator $L$. The problem [5] - [6] has simple boundary conditions and in connection with that is much easier, than the initial problem [1] - [2].
Equality $\lim_{\varepsilon \to 0} \tilde{u}_\varepsilon = \tilde{u}_o$ allows to approve, that the solution of a problem [5] - [6] is as much as close to the solution of a problem [1] - [2] with $\varepsilon \to 0$. This fact allows to name a limiting problem [5] - [6] as regularization of problem [1] - [2].

**Example.** Let there is a homogeneous plate $\Omega$ that has periodic apertures $T_k^\varepsilon$ with radiuses $a(\varepsilon) = \exp\left(-\frac{C_o}{\varepsilon^2}\right)$ where $C_o = const$ (figure 2). Let on borders of a plate and apertures some temperature $T$ is supported, and on its all surface thermal sources operate with density $f$.

The problem is to define temperature $U$ of a plate after the expiration of a significant time interval, i.e. when process is steady. The mathematical model looks like:

$$\Delta U^\varepsilon = -f \text{ On } \Omega_\varepsilon \qquad [7]$$
$$U^\varepsilon\big|_{\partial\Omega} = T \quad U^\varepsilon\big|_{\partial T_k^\varepsilon} = T \quad \left(k = \overline{1,N}\right) \qquad [8]$$

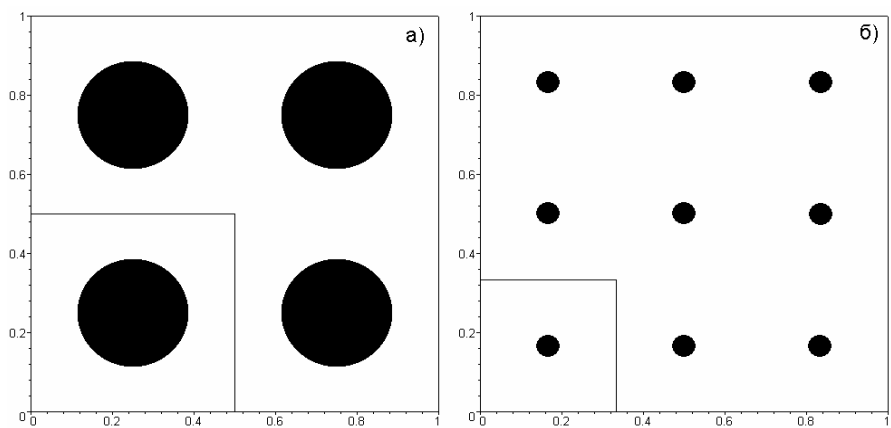

***Fig. 2*** *Location of apertures with radius* $a(\varepsilon) = \exp\left(-\frac{C_o}{\varepsilon^2}\right)$, $C_o = \frac{1}{2}$ *on a plate and length $\varepsilon$ of an internal square, a) $\varepsilon = \frac{1}{2}$ б) $\varepsilon = \frac{1}{3}$*

Using results of works [CM97], [KL01] and the scheme above it is easy to show, that a limiting problem to [7] - [8] with $\varepsilon \to 0$ will be a following problem:

$$-\Delta U + \mu U = f + \mu T \text{ On } \Omega \quad [9]$$

$$U\big|_{\partial\Omega} = T \quad [10]$$

where $\mu = \frac{\pi}{2}\frac{1}{C_o}$.

**The numerical solution of a problem [9] - [10].** For solution of a problem [9] - [10] we will use a method of finite differences.

In Mathcad Professional environment there is a multigrid function which by net method solves a boundary value problem:

$$\Delta G = F \text{ On } \Omega \quad F\big|_{\partial\Omega} = 0$$

I.e. $G = \text{multigrid}(F)$.

Let's consider an implicit iterative scheme:

$$\Delta G_{n+1} + \mu \cdot G_n = -f \quad G_n|_{\partial\Omega} = 0 \quad G_0 = 0, \quad n = 1,2,3,... \quad [11]$$

Using multigrid function we bring implicit iterative process [11] to the obvious scheme:

$$G_{n+1} = \text{multigrid}(\mu \cdot G_n - f) \quad G_0 = 0, \quad n = 1,2,3,... \quad [12]$$

The theorem Iterative process [12] converges, i.e. there is a net function $G$, such, that $\|G_n - G\| \to 0$ with $n \to \infty$ where $\|\ \|$ - means usual euclidean norm.

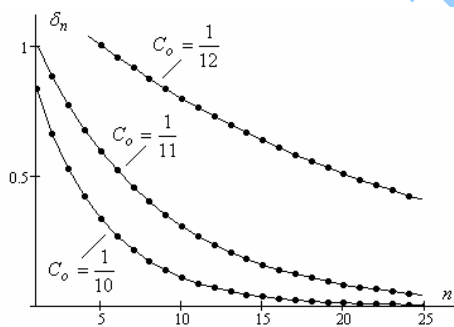

**Fig. 3** *Dependence Graph* $\delta_n = \|G_{n+1} - G_n\|$ *at* $f \equiv 1$ *and* $T = 10$

Considering $U = G + T$, we shall receive, that net function $U$ approximates the solution $U^\varepsilon$ of a boundary value problem [9] - [10].

**Table 1.** *Value of a temperature field* $U$ *in units of a grid with* $16 \times 16$ *size*

$U = $

|    | 0  | 1      | 2      | 3      | 4      | 5      | 6      | 7      | 8      | 9      | 10     | 11     | 12     | 13     | 14     | 15     | 16 |
|----|----|--------|--------|--------|--------|--------|--------|--------|--------|--------|--------|--------|--------|--------|--------|--------|----|
| 0  | 10 | 10     | 10     | 10     | 10     | 10     | 10     | 10     | 10     | 10     | 10     | 10     | 10     | 10     | 10     | 10     | 10 |
| 1  | 10 | 10.004 | 10.006 | 10.007 | 10.008 | 10.008 | 10.008 | 10.008 | 10.008 | 10.008 | 10.008 | 10.008 | 10.008 | 10.007 | 10.006 | 10.004 | 10 |
| 2  | 10 | 10.006 | 10.009 | 10.011 | 10.013 | 10.013 | 10.014 | 10.014 | 10.014 | 10.014 | 10.014 | 10.013 | 10.013 | 10.011 | 10.009 | 10.006 | 10 |
| 3  | 10 | 10.007 | 10.011 | 10.014 | 10.016 | 10.016 | 10.017 | 10.017 | 10.017 | 10.017 | 10.017 | 10.016 | 10.016 | 10.014 | 10.011 | 10.007 | 10 |
| 4  | 10 | 10.008 | 10.013 | 10.016 | 10.017 | 10.018 | 10.019 | 10.019 | 10.019 | 10.019 | 10.019 | 10.018 | 10.017 | 10.016 | 10.013 | 10.008 | 10 |
| 5  | 10 | 10.008 | 10.013 | 10.016 | 10.018 | 10.019 | 10.02  | 10.02  | 10.02  | 10.02  | 10.02  | 10.019 | 10.018 | 10.016 | 10.013 | 10.008 | 10 |
| 6  | 10 | 10.008 | 10.014 | 10.017 | 10.019 | 10.02  | 10.021 | 10.021 | 10.021 | 10.021 | 10.021 | 10.02  | 10.019 | 10.017 | 10.014 | 10.008 | 10 |
| 7  | 10 | 10.008 | 10.014 | 10.017 | 10.019 | 10.02  | 10.021 | 10.021 | 10.021 | 10.021 | 10.021 | 10.02  | 10.019 | 10.017 | 10.014 | 10.008 | 10 |
| 8  | 10 | 10.008 | 10.014 | 10.017 | 10.019 | 10.02  | 10.021 | 10.021 | 10.021 | 10.021 | 10.021 | 10.02  | 10.019 | 10.017 | 10.014 | 10.008 | 10 |
| 9  | 10 | 10.008 | 10.014 | 10.017 | 10.019 | 10.02  | 10.021 | 10.021 | 10.021 | 10.021 | 10.021 | 10.02  | 10.019 | 10.017 | 10.014 | 10.008 | 10 |
| 10 | 10 | 10.008 | 10.014 | 10.017 | 10.019 | 10.02  | 10.021 | 10.021 | 10.021 | 10.021 | 10.021 | 10.02  | 10.019 | 10.017 | 10.014 | 10.008 | 10 |
| 11 | 10 | 10.008 | 10.013 | 10.016 | 10.018 | 10.019 | 10.02  | 10.02  | 10.02  | 10.02  | 10.02  | 10.019 | 10.018 | 10.016 | 10.013 | 10.008 | 10 |
| 12 | 10 | 10.008 | 10.013 | 10.016 | 10.017 | 10.018 | 10.019 | 10.019 | 10.019 | 10.019 | 10.019 | 10.018 | 10.017 | 10.016 | 10.013 | 10.008 | 10 |
| 13 | 10 | 10.007 | 10.011 | 10.014 | 10.016 | 10.016 | 10.017 | 10.017 | 10.017 | 10.017 | 10.017 | 10.016 | 10.016 | 10.014 | 10.011 | 10.007 | 10 |
| 14 | 10 | 10.006 | 10.009 | 10.011 | 10.013 | 10.013 | 10.014 | 10.014 | 10.014 | 10.014 | 10.014 | 10.013 | 10.013 | 10.011 | 10.009 | 10.006 | 10 |
| 15 | 10 | 10.004 | 10.006 | 10.007 | 10.008 | 10.008 | 10.008 | 10.008 | 10.008 | 10.008 | 10.008 | 10.008 | 10.008 | 10.007 | 10.006 | 10.004 | 10 |
| 16 | 10 | 10     | 10     | 10     | 10     | 10     | 10     | 10     | 10     | 10     | 10     | 10     | 10     | 10     | 10     | 10     | 10 |

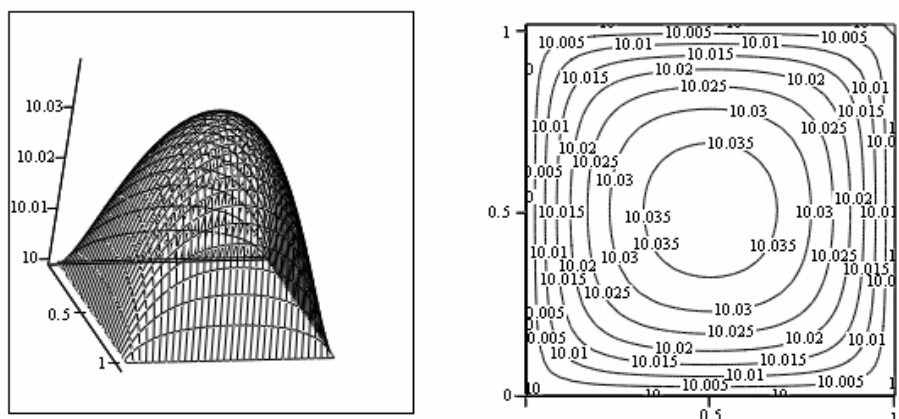

*Fig. 4* *Dependence Graph of a temperature field $U$ and lines of a level for an initial problem at $f \equiv 1$, and $T = 10$*